\documentclass
[11pt]
{article}
\usepackage{latexsym,amssymb,amsmath,amsfonts,latexsym,amscd,amsthm}
\begin{document}
\newtheorem{theorem}{Theorem}
\newtheorem{proposition}{Proposition}
\newtheorem{lemma}{Lemma}
\newtheorem{definition}{Definition}
\newtheorem{corollary}{Corollary}
\setlength{\unitlength}{1mm}

\title{The geometry of the space of leaf  closures of a transversely 
almost K\"ahler foliation}

\author{\normalsize Robert A. Wolak}

\date{Thursday, 2 October 2003}

\maketitle

\begin{abstract}
We study the geometry of the leaf closure space of regular and singular
Riemannian foliations. We give conditions which assure that this leaf space 
is a singular symplectic or K\"ahler space.
\end{abstract}

In recent years physicists and mathematicians working on mathematical 
models of physical phenomena have realised that modelling based on 
geometric structures on now classical smooth manifolds is insufficient. In 
some cases more complicated topological spaces appear naturally and one 
would like to develop geometry of such spaces. One of the well-known examples 
is the orbit space of a smooth action of a compact Lie group. Such a space is 
a stratified pseudomanifold of Goresky-MacPherson, cf. \cite{GM,M,D}. This fact
has been used to describe the topology and structure of the reduced space of 
the momentum map in the singular case, cf. \cite{S1}.

The study of the Riemannian geometry of the orbit space of a smooth action 
of a compact Lie group has been initiated in \cite{AKLM}. One should also 
mention K. Richardson's paper, cf. \cite{RI}, in which the author demonstrates 
that any space of orbits of such an action is homeomorphic to the space of 
the closures of leaves of a regular Riemannian foliation, and, obviously, cf. 
\cite{MO-B}, vice versa.

In this paper we innitiate the study of the geometry of the space of leaf 
closures of a Riemannian foliation, both regular or singular. The first part
is concerned with the space of leaf closures of an interesting class of regular 
Riemannian foliations - transversely almost hermitian, which includes 
transversely almost K\"ahler foliations. In this case the foliated manifold is 
presymplectic. The aim of this note is to describe the geometry of the 
leaf space $M/\overline{\mathcal F}$; in other word to describe the way in 
which the  Riemannian, almost complex and symplectic structures descend 
onto this leaf space.

One of the reasons of the study of such foliations is the fact that recently 
there have been a renewed interest in non-integrable geometries 
associated to Riemannian structures - in fact these geometries correspond to 
the choice of an additional geometric structure compatible with the 
Riemannian metric, cf. \cite{FR}. Almost Hermitian  structures are one of the 
best known examples of such structures. Foliated (almost) Hermitian manifolds 
or foliated K\"ahler manifolds are of interest as they combine three foliated 
structures: a Riemannian, an almost complex and a symplectic ones. These 
structures appear quite naturally in geometry as Sasakian manifolds form 
a special class manifolds foliated by transversally K\"ahler 1-dimensional 
Riemannian foliations, cf. \cite{WO-D} and $\mathcal{K}$-manifolds give 
another example of such foliated manifolds, cf. \cite{DKPW}.

\section{Geometric structures on foliated manifolds} 
It is well-known that on a smooth manifold $M$ of dimension $m$

- the existence of a Riemannian metric $g$ is equivalent to the existence
 of an $O(m)$-reduction $B(M,O(m))$ of the frame bundle $L(M)$ of $M$;
 
- the existence of an almost complex structure $J$ is equivalent to the 
existence  of a $GL(n, {\mathbb C})$-reduction $B(M,GL(n, {\mathbb C}))$ 
of the frame bundle $L(M)$ of $M$, $m=2n$;

- the existence of a 2--form $\omega$ of maximal rank is equivalent to the 
existence  of an $Sp(n)$-reduction $B(M,Sp(n))$ of the frame bundle $L(M)$ 
of $M$, $m=2n$;

\medskip

The almost complex structure $J$ is a complex structure iff

- the Nijenhuis tensor $N(J)(X,Y)=[X,Y] + J[X,JY] + J[JX,Y]+ [JX,JY]=0;$

\noindent
or

- the structure $B(M,GL(n, {\mathbb C}))$ is integrable;

\noindent
or

- the structure tensor of $B(M,GL(n, {\mathbb C}))$ vanishes identically.

\medskip

The 2--form $\omega$ is closed iff 

\noindent
- the structure tensor of $B(M,Sp(n))$ vanishes identically;

\noindent
or

\noindent
- the $Sp(n)$-structure $B(M,Sp(n))$ is integrable.

\medskip

Any two of the three groups $GL(n, {\mathbb C}), Sp(n)$ and $O(2n)$ intersect 
giving $U(n)$ - the maximal compact subgroup of $GL(n, {\mathbb C})$ and 
$Sp(n)$.

\medskip
Any $U(n)$--reduction $B(M,U(n))$ of $L(M)$ assures the existence of a Riemannian metric $g$ 
and an almost complex structure $J$ such that

$$g(Ju,Jv) = g(u,v).$$

\noindent
The associated 2--form $\omega (u,v) = g(Ju,v)$ is of maximal rank. 

\medskip

If the structure tensor of $B(M,U(n))$ vanishes identically, then the structure 
tensors of the corresponding $B(M,GL(n, {\mathbb C}))$ and  $B(M,Sp(n))$ 
 vanish identically, thus the associated almost complex structure $J$ is 
 integrable and the associated 2--form $\omega$ is closed. Therefore such a
structure is a K\"ahler structure.

\medskip

In the foliated case we have the following.
Let the foliation ${\cal  F}$ be  given by a cocycle ${\cal U}=\{U_i,f_i, 
g_{ij}\}$  modelled  on  a  manifold $N_0$ of dimension $2q$, i.e.

i) $\{U_i\}$ is an open covering of $M$,

ii) $f_i \colon U_i \longrightarrow N_0$ are submersions with
connected fibres defining $\cal F$,

iii) $g_{ij}$ are local diffeomorphisms of $N_0$ and
$g_{ij}\circ f_j = f_i$ on $U_i \cap U_j$.

    The manifold $N = \coprod f_i(U_i)$ we call  the  transverse  manifold  of
${\cal F}$ associated to the cocycle ${\cal U}$  and  the  pseudogroup  $H$
generated  by  $g_{ij}$  the  holonomy  pseudogroup  (representative)  on  the
transverse manifold $N$.

It is not difficult to verify, cf. \cite{WO-T,WO-S}, that the following 
conditions a,b,c and d are equivalent:

\noindent
a) On the manifold $N$ there exist

i) a positive definite symmetric tensor $\overline{g}$;

ii) an almost complex structure $\overline{J}$  such that:
 
for any $u,v\in TN \;\;\; \overline{g}(\overline{J}u,\overline{J}v) = 
\overline{g}(u,v)$;

for any $h \in {\mathcal H}$ and any  $u,v\in TN \;\;\; 
\overline{g}(dhu,dhv) = \overline{g}(u,v)$;

for any $h \in {\mathcal H}\;\;\; dh\overline{J}= \overline{J}dh$;

\noindent
then  the 2--form $ \overline{\omega }, \;\; \overline{\omega }(u,v) = 
\overline{g}(\overline{J}u,v)$ 
for any $u,v \in TN, $ is $\mathcal H$-invariant and  of maximal rank. 
Therefore $(N,\overline{g},\overline{J})$ is an almost hermitian manifold and 
elements of the holonomy pseudogroup are holomorphic isometries.

\noindent
b) On the normal bundle $N(M,{\mathcal F})$ there exist

i) an almost complex structure $J$,

ii) a positive definite symmetric tensor $g$ such that 

for any $u,v\in N(M,{\mathcal F}) \;\;\; g(Ju,Jv) = g(u,v)$;

for any $X \in T{\mathcal F} \;\;\; L_Xg =0 \;\;\; L_XJ=0$;

\medskip

Between these three structures we have the following relations (on the 
foliated level), cf. \cite{AC} for the non-foliated case.

1) If a base-like 2--form $\omega$ of maximal rank and a foliated 
almost complex structure are given such that 

$$\omega (Ju,Jv) = \omega (u,v)$$
$$\omega(Ju,u) > 0 \;\;  \forall u \in N(M,{\mathcal F}), \; \; u\not= 0,$$

\noindent
then $g(u,v) = \omega(u,Jv)$ is a foliated positive definite Riemannian 
metric on $N(M,{\mathcal F}).$

2) If a foliated Riemannian metric $g$ and a foliated almost complex 
structure $J$ are given on $N(M,{\mathcal F})$ such that $g(Ju,Jv) = g(u,v)$ 
for any $u,v \in N(M,{\mathcal F}),$ then $\omega (u,v) = g(Ju,v)$ is a 
non--degenerate foliated 2--form of maximal rank, perhaps not closed.

3) If a foliated Riemannian metric $g$ and a foliated 2--form $\omega$ of 
maximal rank on $N(M,{\mathcal F})$ are given, then there exists a compatible 
foliated almost complex structure $J$ on $N(M,{\mathcal F}).$ However, in 
general, the Riemannian metric $\tilde{g}(u,v) = \omega(u,Jv)$ is different
from the initial Riemannian metric $g$. To prove this fact we have to mimmick the 
considerations in \cite{AC}, pp.68-69.

\medskip

\noindent
c) Let $L(M,{\mathcal F})$ be the frame bundle of $N(M,{\mathcal F})$.
$L(M,{\mathcal F})$ admits a foliated reduction, cf. \cite{MO-B,WO-T,WO-S}, to
the structure group $U(q)$.

\medskip

\noindent
d) The frame bundle $L(N)$  admits a $U(q)$-reduction which is 
$\mathcal H$-invariant.

\medskip

The vanishing of the normal Nijenhuis tensor   $N(J)(X,Y) = [X,Y] +J[JX,Y] + 
J[X,JY] -[JX,JY] = 0 \; mod {\mathcal F}$ for any $X,Y \in 
N(M,{\mathcal F}),$ ensures that the almost complex structure on 
$N(M,{\mathcal F})$ is a complex one. This condition is equivalent to the 
vanishing of the Nijenhuis tensor $N(\overline{J})$ of the induced almost 
complex structure $\overline{J}$ on the transverse manifold $N$.

 The same is true for the corresponding reduction $B(N,Sp(q))$ and the 
2--form $\overline{\omega}$. The 2-form $\omega(u,v) = g(Ju,v)$ is a  basic 
2--form inducing on the normal bundle $N(M,{\mathcal F})$ a non-degenerate 
2--form.

Let $B(M,Sp(q);{\mathcal F})$ be the $Sp(q)$--reduction of $L(M;{\mathcal F})$ 
corresponding to $\omega.$ The fact that the form is base--like translates 
itself into the condition "foliated reduction". Therefore the condition 
"the structure tensor vanishes identically" is equivalent to "$d\omega = 0$".
The same is true for the corresponding reduction $B(N,Sp(q))$ and the 2--form 
$\overline{\omega}$; the reduction should be holonomy invariant, cf. 
\cite{WO-T,WO-S}.

\section{The regular case: the topology of the leaf closure space 
$M/\overline{\mathcal F}$}

An almost transversely hermitian foliation is, in particular, a Riemannian 
foliation, therefore the manifold $M$ admits the following stratification, 
cf. \cite{MO-B}.

Let k be any number between 0 and n. Define

$$\Sigma_k = \{ x\in M: x \in L_{\alpha}, dimL_{\alpha} = k
\}.$$

\noindent
The leaves of $\mathcal F$ in any $\Sigma_k$ are of the same dimension but 
they 
can have different holonomy. P. Molino demonstrated that connected components 
of these subsets are submenifolds of $M$ and that $\overline{\Sigma}_k \subset 
\bigcup _{i \leq k} \Sigma_i.$ For some i the sets $\Sigma_i$ may be empty.
Let $k_0$ be the maximun dimension of leaves of $\mathcal F$. Then the set 
$\Sigma_{k_0}$ is open and dense in $M$. It is called the principal stratum.

In many respects the foliation ${\mathcal F}\vert \Sigma_i$ behaves as a RF on 
a compact manifold, in particular, the closure of any its leaves is in 
$\Sigma_i.$  There we can define

$$ \Sigma_{ij} =\{ x \in \Sigma_i \colon x \in L \in {\mathcal F}, \; dim
\overline{L} = j\}.$$

Each $\Sigma_{ij}$ is a submanifold of $\Sigma_i.$
The closures of leaves of $\mathcal F$ induce a regular RF ${\mathcal F}_{ij}$ 
of compact leaves on ${\Sigma}_{ij}.$ The holonomy group of each 
$\overline{L}$ 
is finite, and there is a finite number of types of the holonomy groups. 
Therefore each set ${\Sigma}_{ij}$ decomposes into a finite number of 
submanifolds 

$$\Sigma_{pj{\alpha}}= \{ x\in L \in {\mathcal F}\colon dimL = p, \; dim\overline{L} 
= p, h(\overline{L},x)  \in  \alpha\}$$
\noindent
where $h(\overline{L},x) $ is the holonomy group of $\overline{L}$ in 
$\Sigma_{pj}$ and $\alpha$ is conjugacy class of finite subgroups of $O(q_{pj})$
where $q_{pj}$ is the codimension of $\overline{\mathcal F}$ in $\Sigma_{pj}$.
In this way, we have obtained a stratification ${\mathcal S}= \{\Sigma_{\gamma}\}$ 
of $(M,{\mathcal F})$ into submanifolds on which both $\mathcal F$ and 
$ \overline{\mathcal F}$ define regular Riemannian foliations and the foliation
$ \overline{\mathcal F}$ is without holonomy. The stratification introduced above 
is slightly finer than the stratification defined by M. Pierrot in \cite{P}.
On each stratum $S_{\alpha}$ of $\mathcal S$ the foliation 
$\overline{\mathcal F}$ defines a  compact RF without holonomy, so the natural 
projection on the corresponding stratum $\overline{S}_{\alpha}$ of  
$\overline{\mathcal S}$ is a locally trivial fibre bundle. 

Let us choose a stratum $\Sigma_{\alpha} \in {\mathcal S}$. The foliation 
$\overline{\mathcal F}$ is regular on  $\Sigma_{\alpha},$ so we have the 
following orthogonal splitting of the bundle $TM\vert {\Sigma}_{\alpha}:$

$$ (*) \;\;\;\;\; TM = T{\mathcal F} \oplus Q_1 \oplus Q_2 \oplus Q_3$$

\noindent
where $T{\mathcal F} \oplus Q_1 = \overline{\mathcal F}$ and 
$T{\mathcal F} \oplus Q_1 \oplus Q_2 =T\Sigma_{\alpha}.$

\subsection{Symplectic structure}

A foliation $\mathcal F$ is transversally symplectic if it admits a basic 
closed 2-form $\omega$ of maximal rank. Therefore the codimension of  
$\mathcal F$ is even, say $2q$. 
In particular, any transversally almost K\"ahler foliation is 
transversally symplectic. 

{\it If $\mathcal F$ is transversely symplectic, then $C^{\infty}(M, {\mathcal F})
= C^{\infty}(M, \overline{\mathcal F}) = C^{\infty}(M/\overline{\mathcal F})$ 
admits the structure of a Poisson algebra. }

The form $\omega$  projects on the holonomy invariant 
form $\overline{\omega}$, which is 
a symplectic form of $N$.  The symplectic form 
$\overline{\omega}$ defines a Poisson structure $\{,\}_N$ on $N$ which assigns 
to two $H$-invariant functions an $H$-invariant function. Therefore this 
Poisson structure lifts to a Poisson structure $\{,\}_B$ on the set of basic 
functions $C^{\infty}(M,{\mathcal F})$.  The basic functions are the same for both 
foliations $\mathcal F$ and $\overline{\mathcal F}$, therefore $\{,\}_B$ is 
a Poisson structure on the set  $C^{\infty}(M,\overline{\mathcal F})$. This algebra 
of smooth functions can be considered as the smooth structure on the singular 
space $M/\overline{\mathcal F},$ which is compatible with the definition of a 
smooth structure on the orbit space, cf. \cite{S1}.

We would like to prove that our stratified pseudomanifold 
$M/\overline{\mathcal F}$ is a symplectic stratified space, cf. \cite{S1}. 

{\it Let $\mathcal F$ be transversally almost K\"ahler for the Riemannian metric 
$g,$ the almost complex structure $J$ and the symplectic form $\omega$ on the 
normal bundle $N(M,{\mathcal F}).$ If $J(T\overline{\mathcal F}/T{\mathcal F})
\subset T\overline{\mathcal F}/T{\mathcal F},$ then  $M/\overline{\mathcal F}$ 
is a singular symplectic space.}

It is well-known that the closures of leaves of a Riemannian foliations are 
the orbits of the commuting sheaf of this foliation, cf. \cite{MO-B}, which is 
defined using the bundle of transverse orthonormal frames. The local vector 
fields of the commuting sheaf are local Killing vector fields of the induced 
Riemannian metric. In our case we can refine the definition, cf. 
\cite{WO-MA,WO-S}. The compatible foliated Riemannian and symplectic 
structures define a foliated  $U(q)$-reduction $B(M,U(q);{\mathcal F})$ of the 
bundle $L(M,{\mathcal F})$ of transverse frames, i.e. the frames of the normal 
bundle $N({\mathcal F})$. The group $U(q)$ is of type 1, so the foliation 
${\mathcal F}_1$ of the total space of $B(M,U(q);{\mathcal F})$ is transversely 
parallelisable (TP) and the closures of leaves form a regular foliation with 
compact leaves. The projections of these leaves onto $M$ are the closures of 
leaves of $\cal F$. From the general theory of TP foliations, cf. 
\cite{MO-B}, we know that these closures are the orbits,  in the foliated 
sense, of local vector 
fields commuting with  global foliated vector fields, 
in particular with vector fields of the transverse parallelism. These vector 
fields are the lifts to $B(M,U(q);{\mathcal F})$ of local foliated vector 
fields on $(M,{\mathcal F}),$ which are infinitesimal automorphisms of the 
transverse $U(q)$-structure, so they preserve both the transverse Riemannian 
metric, almost complex structure and the associated 2--form.

We shall look at transverse symplectic and holomorphic structures and check 
whether some  of their components project onto the leaf space 
$M/\bar{\mathcal F}$ and verify what structures they induce.

First consider the principal stratum $\Sigma_0$. The splitting of $TM$ reduces
 itself to $T{\mathcal F} \oplus Q_1 \oplus Q_2 $ as $\Sigma_0$ is an open and 
 dense subset $M$. Since $J(Q_1) \subset Q_1,$ then $J(Q_2) \subset Q_2$. 
 Therefore as the splitting is $g$-orthogonal, it is also $\omega$-orthogonal 
 and the transverse symplectic form $\omega$ can be written as 
 $\omega = \omega^{2,0} + \omega^{0,2}$ where $\omega^{2,0}$ and $\omega^{0,2}$
 are homogeneous components with respect to the splitting. Locally, the 
 subbundle  $Q_1$ is spanned by foliated Killing vector fields $X$ such that 
 $L_X\omega =0.$  Their flows preserve the splitting so 
 
 $$L_X( \omega^{2,0} + \omega^{0,2} ) = L_X\omega^{2,0} + L_X\omega^{0,2}$$
 
 \noindent
 and
 \noindent
 $L_X\omega^{2,0} = 0$ and $L_X\omega^{0,2} =0.$
 
 The restrictions of both forms to $Q_1$ and $Q_2,$ respectively, are of 
 maximal  rank. Therefore to prove that $\omega^{0,2}$ is a transverse 
 symplectic form for  the foliation $\bar{\mathcal F}$ on $\Sigma_0,$ it is 
 sufficient to demonstrate  that $d{\omega}^{0,2} =0.$
 
 In fact, $0=d\omega = d\omega^{2,0} + d\omega^{0,2}
 = d_1\omega^{2,0} + d_2\omega^{2,0} + \partial \omega^{2,0} + d_1\omega^{0,2} 
 + d_2\omega^{0,2}.$
 
 Thus $d_1\omega^{2,0}= 0, \;\;\; d_2\omega^{2,0}=0,\;\;\; \partial\omega^{2,0}
 + d_1\omega^{0,2} = 0, \;\;\; d_2\omega^{0,2}=0.$ Therefore it remains to 
 prove  that $d_1\omega^{0,2}=0.$ Now, for any vector field $X$ of the 
 commuting sheaf   $$L_X\omega^{0,2} = 0 = i_Xd\omega^{0,2} + di_X\omega^{0,2} = 
 i_Xd\omega^{0,2}
 = i_Xd_1\omega^{0,2} + i_Xd_2\omega^{0,2} = i_Xd_1\omega^{0,2}.$$
 
 As these vector fields span the subbundle $Q_1$ we obtain $d_1\omega^{0,2} 
 =0,$  and hence $d\omega^{0,2} =0$. Therefore our foliation 
 $\overline{\mathcal F}$ 
 is transversely symplectic for the 2-form $\omega^{0,2}= \bar{\omega}.$
Thus the projection $\pi_0 \colon \Sigma_0 \rightarrow \overline{\Sigma}_0$
projects the 2-form $\bar{\omega}$ to a symplectic form on $\overline{\Sigma}_0,$
which we denote by the same letter.

Let $\Sigma_{\alpha}$ be any stratum of $(M,{\mathcal F}).$ $\overline{\mathcal F}$ 
induces a regular foliation of no holonomy. 

In \cite{WO-BA} we have proved that global i.a. of $\mathcal F$ are tangent to 
the strata and that the module ${\mathcal X}(M,{\mathcal F})$ of these global 
vector fields is transverse to $\overline{\mathcal F}$ in each stratum. If $X$ 
is an i.a., so is $JX.$ Therefore each stratum is $J$-invariant, i.e.

$J(Q_1\oplus Q_2) \subset Q_1\oplus Q_2,$

\noindent
hence the splitting (*) over any stratum is also $\omega$-orthogonal. In this 
case the standard reasoning ensures that the 2-form $\omega_{\alpha} = 
i_{\alpha}^*\omega,$ where  $i_{\alpha}$ is the inclusion of the stratum 
$\Sigma_{\alpha}$ into $M,$ is a transverse symplectic form of the foliated 
manifold $(\Sigma_{\alpha}, {\mathcal F}).$ Having proved that the same 
considerations as for the principal stratum demonstrate that each stratum of the 
stratification $\overline{\mathcal S}$ is a symplectic manifold.
 

Let $\pi_{\alpha} \colon \Sigma_{\alpha} \rightarrow 
\overline{\Sigma}_{\alpha}$ be the local trivial fibre bundle defining the 
foliation $\overline{\mathcal F}$ on $\Sigma_{\alpha}.$ Clearly, the mapping 
$\pi_{\alpha}$ is a morphism of the Poisson algebras $C^{\infty}({\Sigma}_{ 
\alpha},\overline{\mathcal F}), \{ , \}_{\alpha}$ and 
$C^{\infty}(\overline{\Sigma}_{\alpha}),\overline{\{ , \}}_{\alpha},$ where
the Poisson brackets $\{ , \}_{\alpha}$ and $\overline{\{ , \}}_{\alpha}$ are
defined by $\omega_{\alpha}$ and $\bar{\omega}_{\alpha},$ respectively. 
 
To complete the proof that $M/\overline{\mathcal F}$ is a singular symplectic 
space we have to show that for any $\overline{\Sigma}_{\alpha} \in 
\overline{\mathcal S}$ the inclusion $i_{\alpha} \colon \overline{\Sigma}_{\alpha} 
\rightarrow M/\overline{\mathcal F}$ is a Poisson morphism, i.e.

$$\forall f,g \in C^{\infty}(M/\overline{\mathcal F}) \;\;\;\; \{f,g\}_B\vert 
\overline{\Sigma}_{\alpha} = \overline{\{f\vert \overline{\Sigma}_{\alpha} , 
g \vert \overline{\Sigma}_{\alpha} \}}_{\alpha}.$$

\noindent
In fact,  

$$\{f,g\}_B\vert \overline{\Sigma}_{\alpha} = \omega (X_f,X_g)\vert 
\Sigma_{\alpha}
=\omega_{\alpha}(X_f,X_g) = \omega_{\alpha}(X_f\vert \Sigma_\alpha, X_g\vert 
\Sigma_{\alpha}) = \{ f\vert \overline{\Sigma}_{\alpha}, g \vert 
\overline{\Sigma}_{\alpha} \}_{\alpha}$$

\noindent
as vector fields $X_f, X_g$ are tangent to strata.
 
 \medskip
 The almost complex structure $J_2$ on $Q_2$
 
 $$g(J_2v,w) = \omega(v,w) $$
 
 \noindent
 for any $v,w \in Q_2$ is the restriction of the almost complex structure $J$ 
 to $Q_2$ and therefore its Nijenhuis tensor $N_{J_2} $ is equal to 0.

\section{Singular transversely almost K\"ahler foliations}

The notion of a singular Riemannian foliation was introduced by Pierre Molino 
in \cite{MO-AM}, see also \cite{MO-B}. However, we are not aware
of any notion of an almost complex structure adapted to a singular foliation. 

\begin{definition}

A (1,1)-tensor field $J\in Hom(TM,TM)$ is a foliated almost complex structure
iff

i) $J(T{\mathcal F}) \subset T{\mathcal F};$

ii) for any i.a $X$ of $\mathcal F$ the vector field $JX$ is also an i.a. 
of ${\mathcal F};$

iii) for any i.a $X$ of $\mathcal F$ $J^2X = -X \;\; mod {\mathcal F}.$

\end{definition}

\noindent
{\bf Properties} 

1) In the regular case any transverse almost complex structure can 
be extented to a foliated almost complex structure.

2) If $J(T\overline{\mathcal F}) \subset T\overline{\mathcal F},$ then the 
strata of $(M,{\mathcal F})$ are $J$-invariant provided that the foliation 
$\mathcal F$ is Riemannian.

It is a simple consequence of (i) of the definition and of the fact that 
global i.a. of $\mathcal F$ are transverse to the closures of leaves in the
strata, cf. \cite{WO-BA}.

3) On any stratum ${\Sigma}_{\alpha}$ the foliations $\mathcal F$ and 
$\overline{\mathcal F}$ are regular. Then the tensor field $J_{\alpha} =
J\vert T{\Sigma}_{\alpha}$ is well defined and induces transverse almost
complex structures for ${\mathcal F}\vert {\Sigma}_{\alpha}$.

A foliated Riemannian metric $g$ and a foliated almost complex structure $J$
are said to be compatible if $g(JX,JY) = g(X,Y)$ for any vectors $X,Y \in TM.$
Then on any stratum ${\Sigma}_{\alpha}$ the induced Riemmanian metric 
$g_{\alpha}$ and the induced foliated almost complex structure $J_{\alpha}$ 
are compatible. Then the 2-form ${\omega}_{\alpha} (X,Y) = g_{\alpha} (J_{\alpha} 
X, Y)$ for $X,Y \in T{\Sigma}_{\alpha}$ is of maximal rank.
Moreover, $J_{\alpha}$ induces a transverse almost complex structure for 
$\overline{\mathcal F}\vert {\Sigma}_{\alpha}.$
\begin{definition}
A singular foliation $\mathcal F$ on $M$ is said to be transversely almost 
K\"ahler if it admits a foliated Riemannian metric $g$ and a foliated almost
complex structure, which are compatible, and such that on any stratum 
${\Sigma}_{\alpha}$ of the associated stratification of $M$ the 2-forms 
$\omega_{\alpha}$  are closed. Such a structure is called transversely K\"ahler
if the induced almost complex structures are transversely integrable for 
${\mathcal F}_{\alpha}.$
\end{definition}

\begin{theorem}
Let $\mathcal F$ be a transversally almost K\"ahler singular foliation on a 
compact manifold $M$. Then the space of leaf closures $M/\overline{\mathcal F}$ 
is an almost K\"ahler singular space.
\end{theorem}


\begin{thebibliography}{aa}

\bibitem{AKLM} D. Alekseevsky, A. Kriegl, M. Losik, P. W. Michor, 
The Riemannian geometry of orbit spaces, the metric, geodesics and
integrable systems, Publ. Math. Debrecen 62 (2003), 247-276.
 
\bibitem{BL}  L. Bates, E. Lerman, Proper group actions and symplectic 
stratified spaces, Pacific J. Math. 181 (1997), 201-229 

\bibitem{AC} Ana Cannas da Silva, Symplectic Geometry, Spriger LNM 1764, 2001

\bibitem{D} M. Davis, G-smooth manifolds as collections of fibre bundles, 
 Pacific math J. 77 (1978), 315--363.

\bibitem{DKPW} L. Di Terlizzi, J. Konderak, A.M. Pastore,  R.Wolak,  
K-structures and foliations, Ann. Univ. Sci. Budapest 44 (2001), 171--182.

\bibitem{FR} Th. Friedrich, On types of non-integrable geometries, 
preprint arXiv:math.DG/0205149

\bibitem{GM} M. Goresky, R. MacPherson, Intersection homology theory, Topology 19 (1980), 135-165.

\bibitem{HU} J. Huebschmann, K\"ahler quantizaton and reduction,
arXivmathSG/0207166

\bibitem{M} P. Michor, Isometric Actions of Lie Groups and Invariants, 
 Lecture notes vienna 1997.

\bibitem{MO-AM} P. Molino, D\'esingularisation des feuilletages riemanniens, 
American J. Math. 106 (1984), 1091-1106.

\bibitem{MO-B}  P.  Molino,  Riemannian  Foliations,  Progress  in  Math.  73, 
Birkh\"auser (1988)

\bibitem{P} M. Pierrot, Orbites des champs feuillet\'es pour un feuilletages 
riemanniens sur une vari\'et\'e compacte, C. R. Acad. Sc. Paris 301 (1985),
443-445.

\bibitem{RI} K. Richardson, The transverse geometry of G-manifolds and Riemannian foliations, Illinois J. Math. 45 (2001), 517--535.

\bibitem{SW}  M. Saralegi, R. Wolak,  The BIC for conical fibrations, 
ArXiv.mathDG/0202013.

\bibitem{S1} R.Sjamaar, E. Lerman, Stratified symplectic spaces and reduction, 
Ann. Math. 134 (1991), 375-422.

\bibitem{WO-T}  R. Wolak, Foliated and associated geometric structures on
foliated manifolds, Ann. FAc. Sc. Toulouse 10,3 (1989), 337-360.

\bibitem{WO-S}  R. Wolak, Geometric Structures on Foliated MAnifolds, Santiago de Compostela 1989.

\bibitem{WO-MA} R. Wolak, Foliated G-structures and Riemannian foliations,
Manuscripta Math. 66 (1989), 45--59.

\bibitem{WO-BA} R. Wolak, Pierrot's theorem for singular Riemannian foliations, Publ. Matem. 38 (1994), 433-439

\bibitem{WO-D} R. Wolak, Contact CR-submanifolds in sasakian  manifolds - a foliated  
approach, Publ. Math. Debrecen 56, 1-2 (2000), 7-19
\end{thebibliography}
\end{document}